\documentclass[english]{IEEEtran}
\usepackage[T1]{fontenc}
\usepackage[latin9]{inputenc}
\usepackage{verbatim}
\usepackage{float}
\usepackage{mathrsfs}
\usepackage{amsmath}
\usepackage{amsthm}
\usepackage{amssymb}
\usepackage{graphicx}

\makeatletter

\floatstyle{ruled}
\newfloat{algorithm}{tbp}{loa}
\providecommand{\algorithmname}{Algorithm}
\floatname{algorithm}{\protect\algorithmname}

\theoremstyle{plain}
\newtheorem{thm}{\protect\theoremname}
\theoremstyle{definition}
\newtheorem{problem}{\protect\problemname}
\theoremstyle{remark}
\newtheorem{rem}{\protect\remarkname}
\theoremstyle{plain}
\newtheorem{prop}{\protect\propositionname}
\theoremstyle{plain}
\newtheorem{cor}{\protect\corollaryname}
\theoremstyle{plain}
\newtheorem{lem}{\protect\lemmaname}

\makeatother

\usepackage{babel}
\providecommand{\corollaryname}{Corollary}
\providecommand{\lemmaname}{Lemma}
\providecommand{\problemname}{Problem}
\providecommand{\propositionname}{Proposition}
\providecommand{\remarkname}{Remark}
\providecommand{\theoremname}{Theorem}

\begin{document}
\title{Duality-based Dynamical Optimal Transport of Discrete Time Systems}

\author{Dongjun Wu and Anders Rantzer
	\thanks{*This project has received funding from the European Research
		Council under the European Union's Horizon 2020 research
		and innovation programme under grant agreement No 834142; 
		Excellence Center ELLIIT; Wallenberg AI, Autonomous
		Systems and Software Program (WASP).
}
	\thanks{$^{1}$ The authors are with the Department of Automatic
		Control, Lund University, Box 118, SE-221 00 Lund, Sweden. {\tt\small
	firstname.lastname@control.lth.se}.}
}

\maketitle
\begin{abstract}
We study dynamical optimal transport of discrete time systems (dDOT)
with Lagrangian cost. The problem is approached by combining optimal
control and Kantorovich duality theory. Based on the derived solution,
a first order splitting algorithm is proposed for numerical implementation.
While solving partial differential equations is often required in
the continuous time case, a salient feature of our algorithm is that
it avoids equation solving entirely. Furthermore, it is typical to
solve a convex optimization problem at each grid point in continuous
time settings, the discrete case reduces this to a straightforward
maximization. Additionally, the proposed algorithm is highly amenable
to parallelization. For linear systems with Gaussian marginals, we
provide asemi-definite programming formulation based on our theory.
Finally, we validate the approach with a simulation example.
\end{abstract}

\section{Introduction}

Optimal transport has gained wide popularity in systems and control
in the past a few years \cite{chen2021optimal}. The growing interest
is driven by various applications encompassing uncertainty control
\cite{chen2021controlling}, robotic swarm \cite{krishnan2018distributed}
and nonlinear filtering \cite{taghvaei2021optimal}. Depending on
the scenario, different formulations of optimal transport have been
used, which can broadly be categorized into two types, static and
dynamic. The static formulation, including the Monge, Kantorovich
and the dual problems, aims at finding a static mapping which sends
an initial probability distribution to a target distribution. In contrast,
the dynamic formulation answers not only ``which should move to where'',
but also ``how should a movement be dynamically achieved''. Both types
of problems have been extensively studied in the literature, and comprehensive
treatments of these formulations and their interrelations can be found
in \cite{villani2021topics,ambrosio2021lectures,villani2009optimal}.

In systems and control, the natural setting is the dynamic formulation
initiated by Brenier and Benamou in the seminal work \cite{benamou2000computational},
which can be seen as dynamic optimal transport (DOT) over a first
order integrator with quadratic cost function. Since then, this framework
has been extended to accommodate more complex dynamics and cost structures,
such as integrator with Lagrangian cost \cite{bernard2007optimal,ghoussoub2018optimal},
linear systems with quadratic cost \cite{chen2016optimal}, stochastic
systems \cite{chen2021stochastic}, systems with nonholonomic constraints
\cite{agrachev2009optimal,ambrosio2004optimal} and nonlinear control-affine
systems \cite{elamvazhuthi2023dynamical}. 

A bottleneck in continuous dynamic optimal transport (cDOT) is its
lack of scalability to large-scale systems. Typically, cDOT requires
solving partial differential equations (PDEs). For example, in \cite{benamou2000computational}
and \cite{elamvazhuthi2023dynamical}, solving a Poisson equation
with Neumann boundary conditions is a key step; in \cite{caluya2021wasserstein,chen2021stochastic},
Fokker-Planck or Hamilton-Jacobi equations arise. Moreover, these
algorithms often involve small scale optimization problems at each
grid point of the state space, making computations prohibitive as
the dimension of the state space increases \cite{benamou2000computational,papadakis2014proximal,elamvazhuthi2023dynamical}. 

The study of dynamic optimal transport over discrete time systems
(dDOT) is more recent, and it offers notable advantages over its continuous
counterpart. A key benefit of dDOT is its flexibility. For instance,
cDOT often requires the initial distribution to be absolutely continuous
\cite{benamou2000computational,agrachev2009optimal,elamvazhuthi2023dynamical},\footnote{Some of the theories developed in \cite{elamvazhuthi2023dynamical}
are valid for non-absolutely continuous measures, but the numerical
implementation still needs absolute continuity. } while dDOT is capable of handling measures with atoms, in particular,
systems on discrete domains \cite{solomon2018optimal}. On the other
hand, dDOT can also be leveraged to solve cDOT problems by discretizing
the dynamical system beforehand. This discretize-then-optimize technique
has long been employed in optimal control with PDE constraints \cite{betts2005discretize}.
Its main advantage over optimize-then-discretize -- the commonly
used approach for cDOT -- is its flexibility in handling constraints
\cite{Ghobadi2006}. Surprisingly, this technique has not yet been
fully explored in the context of DOT, despite cDOT being inherently
an optimal control problem of a PDE system. In addition to the flexibility
in handling constraints, we show in this paper that the discretize-then-optimize
approach for DOT is computationally more efficient, in that our algorithm
does not involve solving any algebraic equations, which would have
not been expected when there exist PDE constraints. For example in
\cite{papadakis2014proximal}, dDOT involves solving an algebraic
equation which is essentially a Poisson equation.

Among the existing works on dDOT, linear system with quadratic cost
and Gaussian marginals has served as a standard model in most studies.
For example \cite{de2021discrete,ito2023maximum} and the references
therein explores this setting, where the problem is known to be closely related
to covariance control or control of uncertainty -- a field with a
much longer history \cite{hotz1987covariance}. Recent advancements
in this area can be found in works such as \cite{bakolas2018finite,liu2022optimal}.
Nonlinear dDOT were initiated in \cite{biswal2020stabilization,elamvazhuthi2018optimal,hoshino2020finite}.
In \cite{biswal2020stabilization,elamvazhuthi2018optimal}, the authors
recast the problem as a controlled Markovian processes, solving it
with linear \cite{elamvazhuthi2018optimal} or quadratic programming
\cite{biswal2020stabilization}. In \cite{hoshino2020finite}, the
author replaced the terminal constraint by a soft terminal cost and
solved using normalizing flow techniques. \cite{lavenant2018dynamical}
studied dDOT over discrete surfaces by discretizing cDOT proposed
in \cite{benamou2000computational}. Recently, Terpin et. al. studied 
dDOT using dynamic programming \cite{terpin2024dynamic}. However, duality
was not discussed.

In this paper, we introduce a new approach to dDOT, based on optimal
control and Kantorovich duality theory%
. The contribution of this work is threefold:
\begin{enumerate}
\item \textbf{New Formulations of dDOT}: We propose novel equivalent forms
of dDOT. Notably, we introduce a formulation based on Kantorovich
duality and the Bellman equation, which echoes the continuous time
framework developed in \cite{bernard2007optimal} using Mather theory.
However, our formulation extends to nonlinear systems. This contribution
also unveils the connections between optimal control and optimal transport
in the discrete time setting. 
\item \textbf{Efficient First Order Algorithm}: %
{} Building upon our theoretical framework, we propose a new proximal
splitting algorithm to solve dDOT. In each iteration, explicit solutions
exist, eliminating the need for solving PDE or small scale optimization
problem at each grid point. This is in sharp contrast to the continuous
time case. Moreover, all steps of the algorithm can be trivially parallelized. 
\item \textbf{Application to Linear Systems}: We apply the proposed theory
to solve dDOT for linear systems with quadratic costs and Gaussian
marginals. The result is a semi-definite programming (SDP), providing
new perspectives for this class of problems.
\end{enumerate}
\emph{Notations and terminologies:} For given integers $i\le j$,
denote $\mathcal{I}[i,j]$ as the set $\{i,i+1,\cdots,j\}$. Let $T:X\to Y$
be a measurable mapping and $\mu$ a measure on $X$, then the pushforward
measure of $\mu$ by $T$ is defined as $T_{\#}\mu(A):=\mu(T^{-1}A)$
for all measurable sets $A$ in $Y$. A \emph{Polish space} is a separable
completely metrizable topological space. Given a Polish space $X$,
denote $\mathscr{P}(X)$ the space of probability distributions on
$X$.

\section{Theoretical foundations}

Let $X$ be a Polish space and $U$ a compact measurable space. A
discrete time dynamical system on $(X,U)$ is defined by a sequence
of mappings $f_{k}:X\times U\to X$: 
\begin{equation}
x_{k+1}=f_{k}(x_{k},u_{k})\label{eq:sys}
\end{equation}
where $x_{k}\in X,$ $u_{k}\in U$ for all $k\in\mathbb{N}$. In this
paper, we shall study finite horizon problems. Let $T$ be a terminal
time stamp, then the state is defined on for $k\in\mathcal{I}[1,T]$
and the input is defined for $k\in\mathcal{I}[1,T-1]$. We call $\{u_{k}\}_{k=1}^{T-1}$
the input sequence and $\{x_{k}\}_{k=1}^{T}$ the state trajectory. 

Given $x,y\in X$, an optimal control problem on $\mathcal{I}[i,j]$
of Lagrangian form is defined as 
\begin{equation}
c_{i}^{j}(x,y):=\min_{u}\sum_{k=i}^{j-1}L_{k}(x_{k},u_{k})\label{pb:oc}
\end{equation}
for some $i\le j\in\mathbb{N}$ such that $x_{i}=x$, $x_{j}=y$,
where $L_{k}:X\times U\to\mathbb{R}_{\ge0}$ is the associated (non-negative)
Lagrangian. When there is no control input $\{u_{k}\}_{k=i}^{j-1}$
satisfying $x_{i}=x$ and $x_{j}=y$, define $c_{i}^{j}(x,y)=+\infty$.
We assume %
the following standing assumption:

\textbf{Assumption H1}: The mappings $f_{k}$ and $L_{k}$ are continuous
on $X\times U$. The mapping $(x,y)\mapsto c_{1}^{T}(x,y)$ is continuous
and finite for all $(x,y)\in X\times X$. Furthermore, assume that
each minimizer $\{u_{k}\}_{1}^{T-1}$, as a function of $(x,y)$,
is measurable. 

Let $\mu_{1}$, $\mu_{T}$ be two probability measures on $X$, standing
for the initial and target distributions respectively. Assume:

\textbf{Assumption H2}: There exist two functions $(g,h)\in L^{1}(\mu_{1})\times L^{1}(\mu_{T})$
such that $c_{1}^{T}(x,y)\le g(x)+h(y)$ for all $(x,y)\in X\times X$.

{} We consider the following three optimal transport problems. 
\begin{problem}[Monge Problem]
\label{pb:M}
\[
\mathcal{M}(\mu_{1},\mu_{T}):=\inf_{S\text{ measurable}}\int_{X}c_{1}^{T}(x,S(x)){\rm d}\mu_{1}(x)
\]
such that $S_{\#}\mu_{1}=\mu_{T}$. When there exists no $S$ such
that $S_{\#}\mu_{1}=\mu_{T}$, define $\mathcal{M}(\mu_{1},\mu_{T})=+\infty$.
\end{problem}
\begin{problem}[Kantorovich Problem]
\label{pb:K}
\[
\mathcal{K}(\mu_{1},\mu_{T})=\inf_{\pi\in\Gamma(\mu_{1},\mu_{T})}\int_{X\times X}c_{1}^{T}(x,y){\rm d}\pi(x,y)
\]
where $\Gamma(\mu_{1},\mu_{T})$ is the set of probability measures
on \textbf{$X\times X$} whose two marginal distributions are $\mu_{1}$
and $\mu_{T}$ respectively. 
\end{problem}
\begin{problem}[Kantorovich-Rubinstein dual]
\label{pb:KR-dual}
\[
\mathcal{KR}(\mu_{1},\mu_{T})=\sup_{\phi,\varphi}\int_{X}\phi(x){\rm d}\mu_{1}(x)-\int_{X}\varphi(y){\rm d}\mu_{T}(y)
\]
where $\phi\in L^{1}(\mu_{1})$, $\varphi\in L^{1}(\mu_{T})$ and
$\phi(x)-\varphi(y)\le c_{1}^{T}(x,y)$. 
\end{problem}
It is standard results that Problems \ref{pb:M} to \ref{pb:KR-dual}
are equivalent in the sense that they give the same optimal cost \cite{villani2021topics}.
However, all three problems rely on knowing \emph{a prior }the cost
$c_{1}^{T}(x,y)$ for every $(x,y)\in X^{2}$. In other words, one
has to solve the optimal control problem (\ref{pb:oc}) for each pair
$(x,y)\in X$ in the first place, which, except in very limited cases,
e.g., linear system with quadratic Lagrangian, is infeasible in general.
To circumvent this issue, we derive equivalent dynamic formulations
of \ref{pb:M} to \ref{pb:KR-dual}. Similar to cDOT proposed \cite{benamou2000computational},
such formulations will allow us to derive more tractable algorithms
to compute the optimizer.

Define $\mathcal{X}_{U}:=\{(x_{1},u_{1},\cdots,u_{T-1}):x_{1}\in X,\;u_{i}\in U\}$.
An element in $\mathcal{X}_{U}$ is called a curve. Given a curve
$\gamma\in\mathcal{X}_{U}$, the action $\mathcal{A}(\gamma)$ of
$\gamma$, is defined as 
\[
\mathcal{A}(\gamma)=\sum_{i=1}^{T-1}L_{k}(x_{i},u_{i}).
\]
Define two boundary mappings $e_{1},e_{T}:\mathcal{X}_{U}\to X$ as
$e_{1}\gamma=x_{1}$ and $e_{T}\gamma=x_{T}$ where $x_{T}$ is the
terminal state under the control action $\{u_{k}\}_{1}^{T-1}$. We
say that $\gamma$ is minimizing if $\mathcal{A}(\gamma)\le\mathcal{A}(\gamma')$
for all $\gamma'\in\mathcal{X}_{U}$ with $e_{1}\gamma'=e_{1}\gamma$
and $e_{T}\gamma'=e_{1}\gamma$. The set of minimizing curves is denoted
as $\mathcal{GX}_{U}$. 
\begin{rem}
The set $\mathcal{GX}_{U}$ is a reminiscent of space of geodesics
in metric space. 
\end{rem}
Let $\mathscr{P}(\mathcal{X}_{U})$ be the space of probability measures
on $\mathcal{X}_{U}$. Then we are ready to state our first version
of dDOT, which we call lifted dDOT, as the optimization variable is
a probability distribution on the space of curves.
\begin{problem}[Lifted dDOT]
\label{pb:D-curvspace}Consider the minimization problem
\begin{equation}
\mathcal{D}_{1}(\mu_{1},\mu_{T}):=\inf_{\eta\in\mathscr{P}(\mathcal{X}_{U})}\int_{\mathcal{X}_{U}}\mathcal{A}(\gamma){\rm d}\eta(\gamma)\label{eq:pb:lift_ot}
\end{equation}
where $\eta$ is such that $(e_{1})_{\#}\eta=\mu_{1}$ and $(e_{T})_{\#}\eta=\mu_{T}$.
\end{problem}
\begin{rem}
Problem (\ref{eq:pb:lift_ot}) can be seen as the discretized counterpart
of 
\[
\inf_{\eta\in\mathscr{P}(C([0,1],X))}\int_{C([0,1],X)}\mathcal{A}'(\gamma){\rm d}\eta(\gamma)
\]
in which $\mathcal{A}'(\gamma)=\int_{0}^{1}L(t,x(t),u(t)){\rm d}t$.
The difference is that the space $\mathscr{P}(C([0,1],X))$ is infinite
dimensional while $\mathscr{P}(\mathcal{X}_{U})$ is finite dimensional. 
\end{rem}

Next, we state another version of dDOT, which is an analog of Benamou
and Brenier's formulation of cDOT \cite{benamou2000computational}.
To that end, we associated a given control sequence $\{u_{k}\}_{k=1}^{T-1}$
with a mapping sequence $S_{k}(x)=x_{k}$ and a sequence of probability
measure $\mu_{k}=(S_{k})_{\#}\mu_{1}$. And we call $\{\mu_{k}\}_{k=1}^{T}$
an interpolation between $\mu_{1}$ and $\mu_{T}$ under the control
action$\{u_{k}\}_{k=1}^{T-1}$, or simply an interpolation when clear
from the context. Consider the following problem:
\begin{problem}[Dynamic OT]
\label{pb:D}%
\[
\mathcal{D}_{2}(\mu_{1},\mu_{T})=\inf_{u_{\cdot},\mu_{\cdot}}\sum_{k=1}^{T-1}\int_{X}L_{k}(x,u_{k}(x)){\rm d}\mu_{k}(x)
\]
where $\{\mu_{k}\}_{k=1}^{T}$ is an interpolation between $\mu_{1}$,
$\mu_{T}$ under the controller $\{u_{k}\}_{k=1}^{T-1}$. 
\end{problem}
The key difference between the above formulation of dDOT and cDOT
is that the interpolation $\mu_{t}$ in cDOT is known to solve the
continuity equation
\[
\partial_{t}\mu_{t}+{\rm div}(\mu_{t}(f(x,u))=0
\]
which is linear in $\mu_{t}$ (for $u$ fixed), and has been key in
deriving cDOT algorithms. However, no such equation is available in
the discrete time setting. Instead, this PDE is replaced by the relation
$(f_k(x,u_k))_\# \mu_k = \mu_{k+1}$. 
If we introduce a Lagrangian multiplier $\{\varphi_k\}_1^{T-1}$ to accout this relation, i.e., to solve
\begin{align*}
\inf_{u_{\cdot},\mu_{\cdot}}\sup_{\varphi_{\cdot}} & \sum_{k=1}^{T-1}\int_{X}L_{k}(x,u_{k}(x)){\rm d}\mu_{k}(x)\\
 & \quad+\int_X \varphi_{k}(f_{k}(x,u_{k})){\rm d}\mu_{k+1}-\varphi_{k}(x){\rm d}\mu_{k}
\end{align*}
the problem is still highly nonlinear due to the first integrand on the second line of the above equation.
Unlike in the continuous
case, where one can introduce a momentum varialbe $m = \rho u$ to turn the the problem into a convex
optimization in the varialbe $(\rho ,m)$ (at least for affine-nonlinear system, see
\cite{benamou2000computational,elamvazhuthi2018optimal}), it is not clear how to employ such a technique
in the discrete time setting.
In view of this difficulty, we propose a dual formulation, which is
the most important formulation in this paper. 
\begin{problem}[Dual formulation]
\label{pb:dual}Define
\begin{equation}
\mathcal{C}_{1}(\mu_{1},\mu_{T})=\sup_{v_{\cdot}}\int_{X}v_{1}{\rm d}\mu_{1}-v_{T}{\rm d}\mu_{T}\label{eq:C-dual}
\end{equation}
where $\{v_{k}\}_{k=1}^{T}$ is a continuous solution to the boundary-free
Bellman equation
\begin{equation}
v_{k}(x)=\inf_{u\in U}\{L_{k}(x,u)+v_{k+1}(f_{k}(x,u))\};\label{eq:Bellman}
\end{equation}
and 
\[
\mathcal{C}_{2}(\mu_{1},\mu_{T})=\sup_{v_{\cdot}}\int_{X}v_{1}{\rm d}\mu_{1}-v_{T}{\rm d}\mu_{T}
\]
where $\{v_{k}\}_{k=1}^{T}$ is a continuous sub-solution the boundary-free
Bellman equation, i.e.,
\begin{equation}
v_{k}(x)\le\inf_{u\in U}\{L_{k}(x,u)+v_{k+1}(f_{k}(x,u))\}.\label{ineq:Bellman}
\end{equation}
\end{problem}
\begin{rem}
The constraint (\ref{ineq:Bellman}) is convex (in $v_{\cdot}$) while
(\ref{eq:Bellman}) is not. As we shall see in the following, $\mathcal{C}_{1}=\mathcal{C}_{2}$,
that is the reason why to put them in the same problem formulation. 
\end{rem}
The following is the main theoretical result of this paper, which says that Problem \ref{pb:M} - Problem
\ref{pb:dual} are equivalent in certain sense.
\begin{thm}
\label{thm:D=00003DK=00003DM v1}Let $\mu_{1}$ and $\mu_{T}$ be
two probability measures. Then under Assumption H1-H2: \footnote{Here we have omitted the arguments $(\mu_{1},\mu_{T})$ for ease of
notation.}

1) The inequalities hold
\begin{equation}
\mathcal{M}\ge\mathcal{D}_{2}\ge\mathcal{C}_{2}\ge\mathcal{C}_{1}\ge\mathcal{K}=\mathcal{KR}=\mathcal{D}_{1}.\label{eq:thm:ineq}
\end{equation}
If $\mu_{1}$ is absolutely continuous, then all the inequalities
become equalities.

2)%
{} Optimizers exist for Problems \ref{pb:K}, \ref{pb:KR-dual} and
\ref{pb:dual}. In particular, if $\{v_{k}\}_{k=1}^{T}$ is a maximizer
of Problem \ref{pb:dual}, then $(v_{1},v_{T})$ is a maximizer of
Problem \ref{pb:KR-dual}. Conversely, if $(v_{1},v_{T})$ is an optimizer
of Problem \ref{pb:KR-dual}, then $\{w_{k}\}_{k=1}^{T}$ is an optimizer
of Problem \ref{pb:dual} by solving the Bellman equation (\ref{eq:Bellman})
with boundary conditions $w_{1}=v_{1}$ and $w_{T}=v_{T}$.

3) If $\mu_{1}$ is absolutely continuous and $S$ a minimizer of
Problem \ref{pb:M}, then $({\rm Id}\times S)_{\#}\mu_{1}$ is a minimizer
of Problem \ref{pb:K}. Let $\{u_{k}(x)\}_{k=1}^{T-1}$ be the optimal
controller for problem (\ref{pb:oc}) with $x_{1}=x$, $x_{T}=S(x)$,
then $(\{u_{k}\},\{(S_{k})_{\#}\mu_{1}\})_{k=1}^{T-1}$ is a minimizer
of Problem \ref{pb:D}, in which $S_{k}(x)=f_{k}(x,u_{k}(x))$.

4) $\eta$ is a minimizer of Problem \ref{pb:D-curvspace} if and
only if $\eta$ is supported on $\mathcal{GX}_{U}$, and that if $\eta$
is a minimizer, then $\pi:=(e_{1},e_{T})_{\#}\eta$ is a minimizer
of Problem \ref{pb:K}. 

5) %
The optimal controller is given by
\begin{equation} \label{eq:controller}
u_{k}(x)=\arg\min_{u}\{L_{k}(x,u)+v_{k+1}(f_{k}(x,u))\}
\end{equation}
in which $\{v_{k}\}_{k=1}^{T}$ is an optimizer of Problem \ref{pb:dual}. 
\end{thm}
\begin{IEEEproof}
\textbf{Item 1).} $\mathcal{K}=\mathcal{KR}$ is well-known. We show
\[
\mathcal{M}\ge\mathcal{D}_{2}\ge\mathcal{C}_{2}\ge\mathcal{C}_{1}\ge\mathcal{K}=\mathcal{D}_{1}.
\]

\emph{Step 1).} $\mathcal{M}\ge\mathcal{D}_{2}$. If there exists
no $S$ such that $S_{\#}\mu_{1}=\mu_{T}$, then $\mathcal{M}=\mathcal{D}_{2}=+\infty$
and there is nothing to prove. Henceforth, we fix a mapping $S$ such
that $S_{\#}\mu_{1}=\mu_{T}$. By assumption, for each pair $(x,S(x))$,
there is a control sequence $\{u_{k}(x)\}_{k=1}^{T-1}$ which steers
$x_{1}=x$ to $x_{T}=S(x)$ such that $\sum_{k=1}^{T-1}L_{k}(x_{k},u_{k})=c_{1}^{T}(x,S(x))$.
Define $S_{k}(x)=f_{k}(x,u_{k}(x))$ and $\mu_{k}=(S_{k})_{\#}\mu_{1}$
for $k=1,\cdots,T-1$. Noting that $S_{T-1}\circ S_{T-2}\circ\cdots\circ S_{1}=S$,
we have $\mu_{T}=S_{\#}\mu_{1}$, i.e., $\{\mu_{k}\}_{k=1}^{T}$ is
an interpolation between $\mu_{1}$ and $\mu_{T}$. Meanwhile, 
\begin{align*}
\sum_{k=1}^{T-1}\int_{X} & L_{k}(x,u_{k}(x)){\rm d}\mu_{k}(x)\\
 & =\sum_{k=1}^{T-1}\int_{X}L_{k}(S_{k}(x),u_{k}(S_{k}(x)){\rm d}\mu_{1}(x)\\
 & =\int_{X}\sum_{k=1}^{T-1}L_{k}(x_{k},u_{k}){\rm d}\mu_{1}(x)\\
 & =\int_{X}c_{1}^{T}(x,S(x)){\rm d}\mu_{1}(x).
\end{align*}
That is, for any $S$, we can find $(u_{\cdot},\mu_{\cdot})$ such
that $\int_{X}c_{1}^{T}(x,S(x)){\rm d}\mu_{1}(x)=\sum_{k=1}^{T-1}\int_{X}L_{k}(x,u_{k}){\rm d}\mu_{k}$.
This implies $\mathcal{M}\ge\mathcal{D}_{2}$.

\emph{Step 2).} $\mathcal{D}_{2}\ge\mathcal{C}_{2}$. Given a control
sequence $\{u_{k}(\cdot)\}_{1}^{T-1}$, and the corresponding interpolation
$\{\mu_{k}\}_{1}^{T}$ between $\mu_{1}$ and $\mu_{T}$. Let $\{v_{k}(\cdot)\}_{1}^{T}$
be an arbitrary sub-solution of the Bellman inequality (\ref{ineq:Bellman}).
Then by definition, 
\[
v_{k}(x)-v_{k+1}(f_{k}(x,u_{k}(x))\le L_{k}(x,u_{k}(x)),
\]
for all $k\in\mathcal{I}[1,T-1]$. It follows that
\begin{align*}
\sum_{k=1}^{T-1}\int_{X} & L_{k}(x,u_{k}(x)){\rm d}\mu_{k}(x)\\
 & \ge\sum_{k=1}^{T-1}\int_{X}[v_{k}(x)-v_{k+1}(f_{k}(x,u_{k}(x))]{\rm d}\mu_{k}(x)\\
 & =\sum_{k=1}^{T-1}\int_{X}v_{k}(x){\rm d}\mu_{k}(x)-v_{k+1}(x){\rm d}\mu_{k+1}(x)\\
 & =\int_{X}v_{1}{\rm d}\mu_{1}-v_{T}{\rm d}\mu_{T}.
\end{align*}
Since $v_{\cdot}$ is arbitrary, we conclude that $\mathcal{D}_{2}\ge\mathcal{C}_{2}$.

\emph{Step 3).} $\mathcal{C}_{2}\ge\mathcal{C}_{1}$. Obvious.

\emph{Step 4).} $\mathcal{C}_{1}\ge\mathcal{K}$. Fix a bounded continuous
function $\phi:X\to\mathbb{R}$ and define 
\[
v_{k}(x)=\inf_{y}c_{k}^{T}(x,y)-\phi(y),\;\forall k\in\mathcal{I}[1,T].
\]
Note that $v_{k}$ is well-defined (i.e., $v_{k}(x)>-\infty$ for
all $x\in X$) due to Assumption \textbf{H1}. We claim that $v_{k}$
is a solution to the Bellman equation (\ref{eq:Bellman}). Indeed,
for all $k\in\mathcal{I}[1,T-1]$,
\begin{align*}
v_{k}(x) & =\inf_{y}\inf_{u}L_{k}(x,u)+c_{k+1}^{T}(f_{k}(x,u),y)-\phi(y)\\
 & =\inf_{u}\{L_{k}(x,u)+v_{k+1}(f_{k}(x,u))\}.
\end{align*}
By construction $v_{T}(x)=-\phi(x)$ and $v_{1}(x)=\inf_{y}c_{1}^{T}(x,y)-\phi(y)$,
i.e., $v_{1}$ is the $c_{1}^{T}$-transform of $-v_{T}$ (see Appendix).
Since $\phi$ (or $v_{T})$ is arbitrary, we immediately conclude
that $\mathcal{C}^{1}\ge\mathcal{K}$ invoking Kantorovich duality
\cite[Theorem 5.10]{villani2009optimal}.

\emph{Step 5).} $\mathcal{K}=\mathcal{D}_{1}$. Given $\eta\in\mathscr{P}(\mathcal{X}_{U})$
such that $(e_{1})_{\#}\eta=\mu_{1}$ and $(e_{T})_{\#}\eta=\mu_{T}$,
we have
\begin{align}
\int_{\mathcal{X}_{U}}\mathcal{A}(\gamma){\rm d}\eta(\gamma) & =\int_{\mathcal{X}_{U}}\sum_{k=1}^{T-1}L_{k}(x_{k},u_{k}){\rm d}\eta(\gamma)\nonumber \\
 & \ge\int_{\mathcal{X}_{U}}c_{1}^{T}(e_{1}\gamma,e_{T}\gamma){\rm d}\eta(\gamma)\label{ineq:action}\\
 & =\int_{X\times Y}c_{1}^{T}(x,y){\rm d}(e_{1},e_{T})_{\#}\eta(x,y)\nonumber 
\end{align}
However, $(e_{1},e_{T})_{\#}\eta\in\Gamma(\mu_{1},\mu_{T})$ by definition
$\eta$. Thus $\mathcal{D}_{1}\ge\mathcal{K}$. It remains to show
the converse. For each pair $(x,y)\in X^{2}$, there exists a minimizing
curve $\gamma$ connecting $x$ to $y$ by assumption. Define a mapping
$\Psi:X^{2}\to\mathcal{GX}_{U}$ as $(x,y)\mapsto(x,u_{1},\cdots,u_{T-1})$
where $(x,u_{1},\cdots,u_{T-1})$ is any minimizing curve. By Assumption
\textbf{H1}, $\Psi$ is measurable. Let $\pi\in\Gamma(\mu_{1},\mu_{T})$
and define $\eta=(\Psi)_{\#}\pi$ which is a probability measure on
$\mathcal{X}_{U}$ concentrated on $\mathcal{GX}_{U}$ by construction.
Now
\begin{align}
\int_{\mathcal{X}_{U}}\mathcal{A}(\gamma){\rm d}\eta(\gamma) & =\int_{\mathcal{GX}_{U}}\mathcal{A}(\gamma){\rm d}\eta(\gamma)\nonumber \\
 & =\int_{X\times X}c_{1}^{T}(e_{1}\gamma,e_{T}\gamma){\rm d}\Psi_{\#}\pi(\gamma)\label{eq:action-K}\\
 & =\int_{X\times X}c_{1}^{T}(x,y){\rm d}\pi(x,y)\nonumber 
\end{align}
Thus $\mathcal{D}^{1}\le\mathcal{K}$ and we are done with the inequality
(\ref{eq:thm:ineq}). If $\mu_{1}$ is absolutely continuous, then
$\mathcal{M}=\mathcal{K}$ invoking Pratelli's theorem \cite[Theorem 2.2]{ambrosio2021lectures}.
The proof of \textbf{Item 1)} is complete. 

\textbf{Item 2)}. It is standard result that the minimizer of Problem
\ref{pb:K} exists if $c_{1}^{T}(\cdot,\cdot)$ is lower semi-continuous
\cite[Theorem 1.3]{villani2021topics}, and that the maximization
of Problem \ref{pb:KR-dual} is attained \cite[Theorem 5.10]{villani2009optimal}
under the Assumption \textbf{H2}. Let $(\phi,\varphi)$ be a maximizer
of Problem \ref{pb:KR-dual}. We construct a maximizer of Problem
\ref{pb:dual}. Define $v_{k}(x):=\inf_{u\in U}\{L_{k}(x,u)+v_{k+1}(f_{k}(x,u))\}$
for $k\in\mathcal{I}[2,T-1]$ with boundary condition $v_{T}=\varphi$.
We need first verify the continuity of $\{v_{k}\}_{k=2}^{T-1}$. Given
$x,y\in X$, let $u_{k}^{*}\in\arg\min_{u}\{L_{k}(x,u)+v_{k+1}(f_{k}(x,u)\}$
and $\tilde{u}_{k}\in\arg\min_{u}\{L_{k}(x,u)+v_{k+1}(f_{k}(x,u))\}$,
then 
\begin{align*}
v_{k}(x)-v_{k}(y) & \ge L_{k}(x,u_{k}^{*})-L_{k}(y,u_{k}^{*})\\
 & \quad+v_{k+1}(f_{k}(x,u_{k}^{*})-v_{k+1}(f_{k}(y,u_{k}^{*})).
\end{align*}
Swapping the role of $x$ and $y$, we get 
\begin{align*}
v_{k}(y)-v_{k}(x) & \ge L_{k}(y,\tilde{u}_{k})-L_{k}(y,\tilde{u}_{k})\\
 & \quad+v_{k+1}(f_{k}(y,\tilde{u}_{k}))-v_{k+1}(f_{k}(x,\tilde{u}_{k})).
\end{align*}
By continuity of $v_{T}$, $f_{k}$ and $L_{k}$, the right hand side
of both previous inequalities goes to zero, implying that $v_{T-1}$
is continuous. By induction, we conclude that $\{v_{k}\}_{k=2}^{T}$
are continuous. Now, if we can verify that $\phi$ (the first component
of the minimizer of Problem \ref{pb:KR-dual}) satisfies $\phi(x)\le\min_{u\in U}\{L_{1}(x,u)+v_{2}(f_{1}(x,u))\}$,
then $\{\phi,v_{2},\cdots,v_{T}\}$ is a maximizer of Problem \ref{pb:dual}.
By definition of $\{v_{k}\}_{k=2}^{T-1}$, we know that $v_{2}(x)=\min_{u_{\cdot}}\sum_{k=2}^{T-1}L_{k}(x_{k},u_{k})+\varphi(x_{T})$
with $x_{2}=x$. By definition, we have 
\begin{align*}
\phi(x) & \le\inf_{y}c_{1}^{T}(x,y)+\varphi(y)\\
 & =\inf_{u\in U}\inf_{y}\{L_{1}(x,u)+c_{2}^{T}(f_{1}(x,u),y)+\varphi(y)\}\\
 & \le\inf_{u\in U}\{L_{1}(x,u)+v_{2}(f_{1}(x,u))\}
\end{align*}
which is the desired inequality. The rest of \textbf{Item 2) }is obvious.

\textbf{Item 3)}. The first fact is well-known, i.e., the optimizer
of Problem \ref{pb:K} is concentrated on the graph of $S$ \cite[Theorem 2.3]{ambrosio2021lectures}.
The rest follows from \emph{Step 1).}

\textbf{Item 4)}. In view of the inequality (\ref{ineq:action}),
$\eta$ is a minimizer if and only if $\sum_{k=1}^{T-1}L_{k}(x_{k},u_{k})=c_{1}^{T}(e_{1}\gamma,e_{1}\gamma)$
for $\eta$-almost every $\gamma$. But this equality holds true if
and only if $\gamma\in\mathcal{G}\mathcal{X}_{U}$. By (\ref{eq:action-K}),
we see that $\pi$ is minimizing if $\eta$ is.

\textbf{Item 5)}. This is obvious from the preceding proof. 
\end{IEEEproof}
\begin{rem}
Consider the special case $x_2=u_1$ with Lagrangian $L(x,u)=\|x-u\|^{2}$.
Then the Bellman equation says
\begin{align*}
v_{1}(x) & =\inf_{u}\{\|x - u\|^{2}+v_{2}(u)\}=\inf_{y}\{\|x-y\|^{2}+v_{2}(y)\}
\end{align*}
The right hand side is nothing but the $c$-transform of $-v_{2}$,
by taking $c(x,y)=\|x-y\|^{2}$. It is well-known that $\mathcal{K}(\mu,\nu)=\sup_{v}\int_{X}v{\rm d}\mu-v^{c}{\rm d}\nu$
for $v\in L^{1}({\rm d}\mu)$. On the other hand, the Bellman inequality
implies $v_{1}(x)-v_{2}(y)\le\|x-y\|^{2},\forall x,y.$ Therefore,
Problem \ref{pb:dual} reduces to Problem \ref{pb:KR-dual}. In summary,
Problem \ref{pb:dual} is an extension of the Kantorovich-Rubinstein
duality to dDOT. 
\end{rem}



\begin{rem}
The dual formulation (Problem \ref{pb:dual}) is expressed as a linear program.
However, two critical challenges remain to be addressed:
\begin{itemize} 
		\item \textbf{Curse of dimensionality}: The current algorithm is limited to handling
				low-dimensional systems, as it requires discretization of the continuous state space, which becomes computationally prohibitive for higher dimensions. 
		\item \textbf{Extract control law from the value function}: Even if the value function is close to be
				optimal, the controller obtained by \eqref{eq:controller} may be far from being optimal.
\end{itemize}

Several approaches could be explored to mitigate these challenges, such as approximate dynamic programming and
reinforcement learning \cite{bertsekas2019reinforcement}. These issues will be addressed in a coming work.
\end{rem}
\section{Linear Systems with Gaussian Marginals}

Consider the linear system 
\begin{equation}
x_{k+1}=Ax_{k}+Bu_{k}\label{eq:sys-LTI}
\end{equation}
where the state $x\in\mathbb{R}^{n}$ and input $u\in\mathbb{R}^{m}$.
We consider a quadratic Lagrangian 
\begin{equation}
L(x,u)=x^{\top}Qx+u^{\top}Ru\label{eq:quad-Lag}
\end{equation}
 with $Q\ge0$ and $R>0$. 
\begin{prop}
\label{prop:Lin-Gauss-quad}Let $T\ge n$. Assume the system (\ref{eq:sys-LTI})
is controllable and $A$ invertible. Be $\mu_{1}$, $\mu_{T}$ two
Gaussian distributions and (\ref{eq:quad-Lag}) the Lagrangian. Then
for Problem \ref{pb:KR-dual}, it is sufficient to take $\phi$ and
$\varphi$ as quadratic functions, i.e., $\phi=x^{\top}P_{\phi}x$,
$\varphi=x^{\top}P_{\varphi}x$; for Problem \ref{pb:D}, it is sufficient
in (\ref{ineq:Bellman}) to restrict $\{v_{k}\}_{1}^{T}$ to quadratic
functions $v_{k}(x)=x^{\top}P_{k}x$ with $P_{k}$ symmetric.
\end{prop}
\begin{IEEEproof}
For linear system (\ref{eq:sys-LTI}) with quadratic cost (\ref{eq:quad-Lag}),
it is well-known that $c_{1}^{T}(x,y)$ is quadratic in $[x^{\top},y^{\top}]^{\top}$
when the system is controllable. Hence by Lemma \ref{lem:K-Gauss},
\[
\mathcal{K}(\mu,\nu)=\max\int_{\mathbb{R}^{n}}\phi(x){\rm d}\mu(x)-\varphi(y){\rm d}\nu(y)
\]
where $\phi$ and $\varphi$ are restricted to quadratic functions
subject to constraints $\phi(x)-\varphi(y)\le c_{1}^{T}(x,y)$. The
remaining part is a consequence of Theorem \ref{thm:D=00003DK=00003DM v1}.
\end{IEEEproof}
\begin{thm}
\label{thm:Gauss}Consider the linear system (\ref{eq:sys-LTI}) with
Lagrangian (\ref{eq:quad-Lag}). Let $\mu_{1}\sim\mathscr{N}(0,\Sigma_{1}),$
$\mu_{T}\sim\mathscr{N}(0,\Sigma_{T})$. Then the optimal value of
Problem \ref{pb:M} - \ref{pb:dual} is equal to the optimal value
of the following problem:
\begin{align}
\sup_{P_{\cdot}} & {\rm Tr}(P_{1}\Sigma_{1})-{\rm Tr}(P_{T}\Sigma_{T})\label{eq:sde_Gaussian}
\end{align}
subject to the constraints
\begin{align}
0 & <R+B^{\top}P_{k+1}B\nonumber \\
0 & \le\begin{bmatrix}R+B^{\top}P_{k+1}B & B^{\top}P_{k+1}A\\
A^{\top}P_{k+1}B & A^{\top}P_{k+1}A-P_{k}+Q
\end{bmatrix}\label{eq:SDP2}
\end{align}
If $P_{\cdot}$ is the optimizer of Problem (\ref{eq:sde_Gaussian}),
then the optimal controller is given by 
\[
u(x)=-(R+B^{\top}P_{k+1}B)^{-1}B^{\top}P_{k+1}Ax.
\]
\end{thm}
\begin{IEEEproof}
We use the formulation of Problem \ref{pb:dual}. According to Proposition
\ref{prop:Lin-Gauss-quad}, it is sufficient to restrict to quadratic
functions $\{v_{k}\}_{1}^{T}$. We make the ansatz $v_{k}(x)=x^{\top}P_{k}x$
for some symmetric matrices $\{P_{k}\}_{1}^{T}$. Then we shall compute
the optimal value of the following problem $\max_{v_{\cdot}}{\rm Tr}(P_{1}\Sigma_{1}-P_{T}\Sigma_{T})$
subject to
\[
x^{\top}P_{k}x\le\min_{u}\{x^{\top}Qx+u^{\top}Ru+(Ax+Bu)^{\top}P_{k+1}(Ax+Bu)\}.
\]
Invoking Lemma \ref{lem:SDP}, this is equivalent to
\begin{equation}
\begin{bmatrix}R+B^{\top}P_{k+1}B & B^{\top}P_{k+1}Ax\\
x^{\top}A^{\top}P_{k+1}B & x^{\top}(A^{\top}P_{k+1}A-P_{k}+Q)x
\end{bmatrix}\ge0.\label{eq:inter_gauss}
\end{equation}
The above holds only if $R+B^{\top}P_{k+1}B>0$. Indeed, since $x$
is arbitrary, (\ref{eq:inter_gauss}) implies ${\rm Im}B^{\top}P_{k+1}A\subseteq{\rm Im}(R+B^{\top}P_{k+1}B)$
or ${\rm \ker}(R+B^{\top}P_{k+1}B)\subseteq\ker A^{\top}P_{k+1}B$.
That is, if $(R+B^{\top}P_{k+1}B)z=0$, then $A^{\top}P_{k+1}Bz=0$.
Since $A$ is invertible by assumption, the latter equality implies
$P_{k+1}Bz=0$, leading to $Rz=0$ in the former equality; thus $z=0$,
and $R+B^{\top}P_{k+1}B$ is invertible, as desired. By Schur complement,
\begin{align*}
x^{\top}[ & A^{\top}P_{k+1}A-P_{k}+Q\\
 & -A^{\top}P_{k+1}B(R+B^{\top}P_{k+1}B)^{-1}B^{\top}P_{k+1}A]x\ge0.
\end{align*}
Combining this with $R+B^{\top}P_{k+1}B>0$ and applying Schur complement
again, we get (\ref{eq:SDP2}). 
\end{IEEEproof}

\begin{cor}
\label{cor:Wass}The Wasserstein distance
\begin{equation}
W_{2}^{2}(\mu,\nu):=\min_{\pi\in\Gamma(\mu,\nu)}\int_{\mathbb{R}^{n}\times\mathbb{R}^{n}}\|x-y\|^{2}{\rm d}\pi(x,y)\label{eq:W-Gauss}
\end{equation}
where $\mu\sim\mathscr{N}(0,\Sigma_{1})$, $\nu\sim\mathscr{N}(0,\Sigma_{2})$
is equal to the optimal value of the following problem
\[
\max{\rm Tr}(P_{1}\Sigma_{1}-P_{2}\Sigma_{2})
\]
subject to 
\[
\begin{bmatrix}I+P_{2} & P_{2}\\
P_{2} & P_{2}-P_{1}
\end{bmatrix}\ge0,\quad I+P_{2}>0.
\]
\end{cor}

\begin{IEEEproof}
Consider the discrete time system $x_{2}=x_{1}+u$ with Lagrangian
$L(x,u)=\|u\|^{2}$. Then $c_{1}^{2}(x,y)=\|x-y\|^{2}$ and the Corollary
follows from Theorem \ref{thm:Gauss}. 
\end{IEEEproof}

\section{Splitting algorithm}

\subsection{General problem}

We propose a proximal splitting algorithm to solve Problem \ref{pb:dual}
under the Bellman inequality constraint (\ref{ineq:Bellman}). In
this section we assume that $X$ and $U$ are subsets of Euclidean
spaces and $\mu_{1}$ and $\mu_{T}$ are absolutely continuous with
respect to the Lebesgue measure. Let $\rho_{1}$ and $\rho_{T}$ be
the corresponding densities. Note that the constraints can be equally
written as
\[
v_{k+1}(y)-v_{k}(f_{k}(x,u))+L_{k}(x,u)\ge0
\]
for all $(x,u)\in X\times U$ and $k\in\mathcal{I}[1,T]$. Let $C(X)$
be the space of continuous functions on $X$ and $\mathcal{M}(X\times U)$
the space of Borel measures on $X\times U$. If $v_{k}\in C(X)$,
$f_{k}\in C(X\times U)$ and $L_{k}(x,u)\in C(X\times U)$, then Problem
\ref{pb:dual} can be stated as{\small{}
\begin{align}
\sup_{v_{\cdot}\in C(X)^{T}} & \int v_{1}\rho_{1}-v_{T}\rho_{T}+\inf_{\lambda_{\cdot}\in\mathcal{M}(X\times U)^{T-1}}\label{eq:opt_CM}\\
 & +\sum_{k=1}^{T-1}\int_{X\times U}(v_{k+1}\circ f-v_{k}\circ\pi_{X}+L_{k}){\rm d}\lambda_{k}+I_{\{\lambda\ge0\}}
\end{align}
where $\pi_{X}:X\times U\to X$ is the projection on $X$. }Let $\mathcal{X}=C(X)^{T}$,
$\mathcal{Y}=C(X\times U)^{T-1}$ and define a linear operator $\mathcal{L}:\mathcal{X}\to\mathcal{Y}$
as 
\[
(\mathcal{L}v)_{k}(x,u)=v_{k}(x)-v_{k+1}(f_{k}(x,u)).
\]
For $g\in\mathcal{Y}$ and $\lambda\in\mathcal{Y}^{*}$ (the topological
dual of $\mathcal{Y}$), the pairing $\left\langle f,\lambda\right\rangle _{\mathcal{Y}\times\mathcal{Y}^{*}}$
is defined as $\sum_{k=1}^{T-1}\int_{X\times U}g(x,u)d\lambda(x,u)$.
Then (\ref{eq:opt_CM}) reads
\begin{equation}
-\inf_{v\in\mathcal{X}}\sup_{\lambda\in\mathcal{Y}^{*}}F(v)+\left\langle \mathcal{L}v,\lambda\right\rangle _{\mathcal{Y}\times\mathcal{Y}^{*}}-G(\lambda)\label{eq:saddle_CM}
\end{equation}
in which $F(v)=\int\rho_{T}v_{T}-\rho_{1}v_{1},\;G(\lambda)=\left\langle L,\lambda\right\rangle _{\mathcal{Y}\times\mathcal{Y}^{*}}+I_{\{\lambda\ge0\}}$
are convex functionals. 

Problem (\ref{eq:saddle_CM}) is a saddle point problem. Unfortunately,
we are not aware of any first order algorithms which solve the problem
efficiently. The difficulty lies in the fact that neither of the spaces
$\mathcal{X}$ or $\mathcal{Y}^{*}$ are reflexive, while standard
splitting algorithms require them to be at least reflexive or even
uniform convex \cite{lopez2012forward}. In what follows we study
a special case of problem (\ref{eq:saddle_CM}), and leave the general
problem (\ref{eq:saddle_CM}) for future research. 

\subsection{A special case\label{subsec:spec-case}}

We first consider the simple dynamics $x_{k+1}=f_{k}(x_{k})+u_{k}$
with time-dependent Lagrangian cost. %
{} Note that (\ref{ineq:Bellman}) can be equivalently written 
\[
v_{k+1}(y)-v_{k}(x)+L_{k}(x,y-f_{k}(x))\ge0,\quad\forall x,y\in X.
\]
To reduce notation, let us call $L_{k}(x,y-f_{k}(x))$ as $L_{k}(x,y)$.
Then Problem \ref{pb:dual} can be stated as {\small{}
\begin{align}
\sup_{v} & \int v_{1}\rho_{1}-v_{T}\rho_{T}+\inf_{\lambda}\label{eq:CP0}\\
 & +\sum_{k=1}^{T-1}\int(v_{k+1}(y)-v_{k}(x)+L_{k}(x,y))\lambda_{k}{\rm d}x{\rm d}y+I_{\{\lambda\ge0\}}.
\end{align}
}Note that the difference between (\ref{eq:saddle_CM}) and (\ref{eq:CP0})
is that in the latter, the measure $d\lambda_{k}(x,u)$ is replaced
by $\lambda_{k}(x,y)dxdy$.%

We need some notations to turn this problem into a more workable form.
Define two Hilbert spaces
\begin{align*}
\mathscr{H} & =L_{2}(\mathcal{I}[1,T]\times X),\;\mathscr{G}=L_{2}(\mathcal{I}[1,T-1]\times X\times X)
\end{align*}
on which the inner products are defined as $\left\langle v,u\right\rangle _{\mathscr{H}}=\sum_{k=1}^{T}\int_{X}v_{k}u_{k}{\rm d}x$
and $\left\langle \lambda,\eta\right\rangle _{\mathscr{G}}=\sum_{k=1}^{T-1}\int_{X\times X}\lambda_{k}\eta_{k}{\rm d}x{\rm d}y$.
Define an operator $K:\mathscr{H}\to\mathscr{G}$ 
\begin{equation}
(Kv)_{k}(x,y)=v_{k}(x)-v_{k+1}(y),\label{op:K}
\end{equation}
which is linear bounded (see appendix), and two convex functionals
$F:\mathscr{H}\to\mathbb{R}$ and $G:\mathscr{G}\to\mathbb{R}$
\[
F(v)=\int\rho_{T}v_{T}-\rho_{1}v_{1},\;G(\lambda)=\left\langle \lambda,L\right\rangle _{\mathscr{G}}+I_{\{\lambda\ge0\}}
\]
where $L:\mathscr{G}\to\mathbb{R}$ stands for $L_{k}(x,y)$. Now,
the optimization problem can be written more compactly as
\[
-\inf_{v}\sup_{\lambda}F(v)+\left\langle Kv,\lambda\right\rangle _{\mathscr{G}}-G(\lambda).
\]
This saddle point problem is standard and can be solved using for
example Chambolle-Pock algorithm, outlined in Algorithm \ref{alg:CP}.
\begin{algorithm}
Given $(v^{(1)},\lambda^{(1)})$, iterates
\begin{enumerate}
\item $v^{(n+1)}={\rm prox}_{F}^{\tau}(v^{(n)}-\tau K^{*}\lambda^{(n)})$;
\item $\lambda^{(n+1)}={\rm prox}_{G}^{\sigma}(\lambda^{(n)}+\sigma K(v^{(n+1)}+\theta(v^{(n+1)}-v^{(n)}))$.
\end{enumerate}
\caption{Chambolle-Pock algorithm\label{alg:CP}}
\end{algorithm}
 In this algorithm, $\tau$, $\sigma$ and $\theta$ are some constant
parameters. The algorithm converges when we set $\theta=1$ and $\tau\sigma\|K\|^{2}<1$,
i.e., $4\tau\sigma\mu_{X}<1$. The operators $K^{*}$, ${\rm prox}_{F}^{\tau}$
and ${\rm prox}_{G}^{\sigma}$ can be readily computed as:%
{} $(K^{*}\eta)_{k}(x)=\int\eta_{k}(x,y)-\eta_{k-1}(y,x)dy,\;k\in\mathcal{I}[1,T]$
by adopting the notations $\eta_{1}=\eta_{T}=1$; the two proximal
operators are 
\[
{\rm prox}_{F}^{\tau}(v)_{k}=\begin{cases}
v_{1}+\tau\rho_{1}, & k=1\\
v_{T}-\tau\rho_{T}, & k=T\\
v_{k}, & k\in\mathcal{I}[2,T-1].
\end{cases}
\]
and ${\rm prox}_{G}^{\sigma}(\lambda)=\max\{0,\lambda-\sigma L\}$
respectively.%
{} That is, these operations can be easily implemented. The algorithm
is outlined in Algorithm \ref{alg:CP-1}.
\begin{algorithm}
Given $(v^{(1)},\lambda^{(1)})$, iterate
\begin{enumerate}
\item {\small{}$H_{k}^{(n)}(x)=\int_{X}\lambda_{k}^{(n)}(x,y){\rm d}y$,
$H_{k}^{(n)}(y)=\int_{X}\lambda_{k}^{(n)}(x,y){\rm d}x$.}{\small\par}
\item $v_{k}^{(n+1)}=v_{k}^{(n)}-\tau(H_{k}^{(n)}-H_{k-1}^{(n)})$ for $k\ne1,T$.
\item $v_{1}^{(n+1)}=v_{1}^{(n)}-\tau(H_{1}^{(n)}-\rho_{1})$, $v_{T}^{(n+1)}=v_{T}^{(n)}+\tau(H_{T-1}^{(n)}-\rho_{T})$.
\item $\xi_{k}^{(n+1)}(x,y)=v_{k}^{(n+1)}(x)+\theta[v_{k}^{(n+1)}(x)-v_{k}^{(n)}(x)]-v_{k+1}^{(n+1)}(y)-\theta[v_{k+1}^{(n+1)}(y)-v_{k+1}^{(n)}(y)]$.
\item $\lambda_{k}^{(n+1)}=\max\{0,\lambda_{k}^{(n)}+\sigma(\xi_{k}^{(n+1)}-L_{k})\}$.%
\end{enumerate}
\caption{Duality-based algorithm\label{alg:CP-1}}
\end{algorithm}

\begin{rem}
The most computational intensive part of Algorithm \ref{alg:CP-1}
is Step 1) which involves two numerical integrations. Step 5) is a
simple maximization which can be done efficiently. The rest of the
steps are simple arithmetic manipulations or substitutions. All the
steps can be parallelized. 
\end{rem}

\section{Simulation example}

\begin{figure}[h]
\begin{centering}
\includegraphics[scale=0.4]{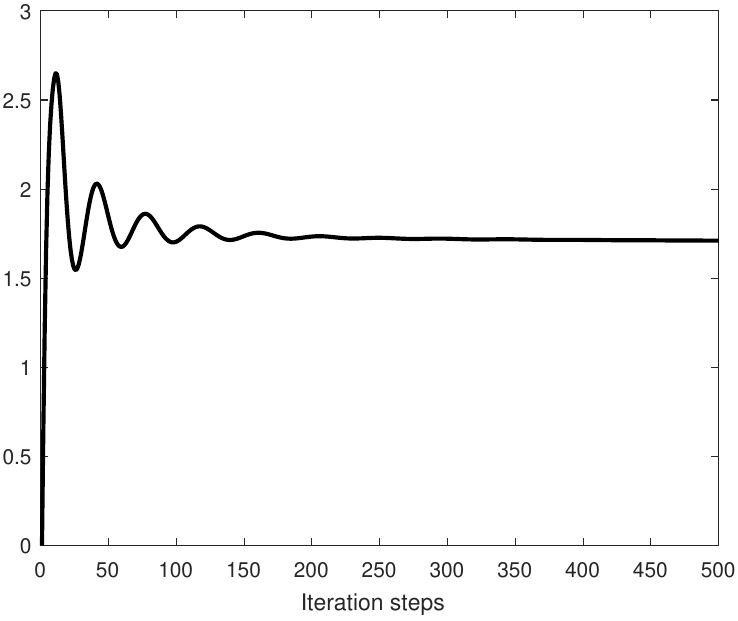}
\par\end{centering}
\caption{The optimal cost.\label{fig:cost}}
\end{figure}

\begin{figure}[h]
\begin{centering}
\includegraphics[scale=0.35]{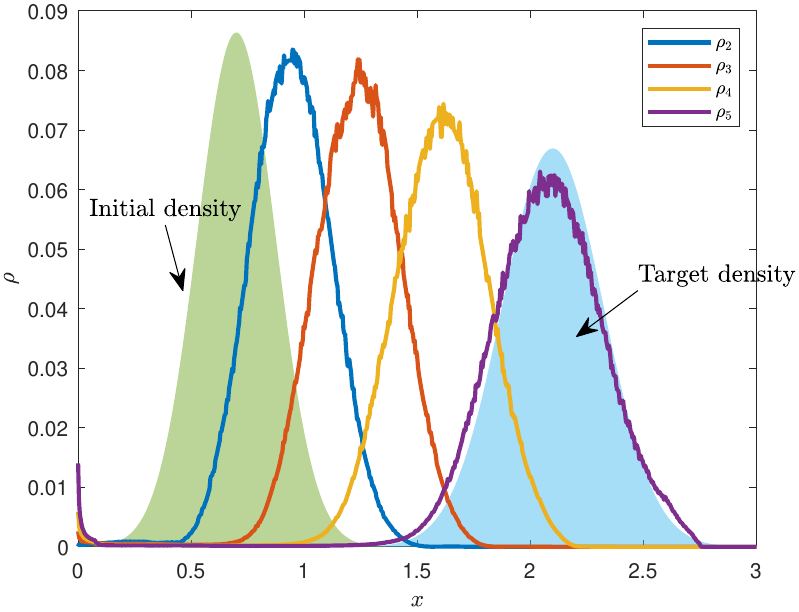}
\par\end{centering}
\caption{The evolution of densities.\label{fig:density} }
\end{figure}

\begin{figure}[h]
\begin{centering}
\includegraphics[scale=0.4]{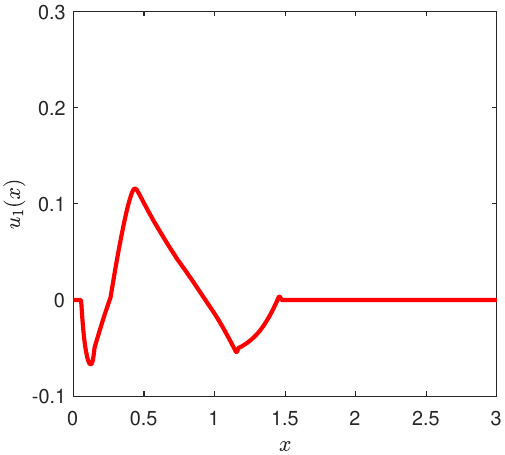}\includegraphics[scale=0.4]{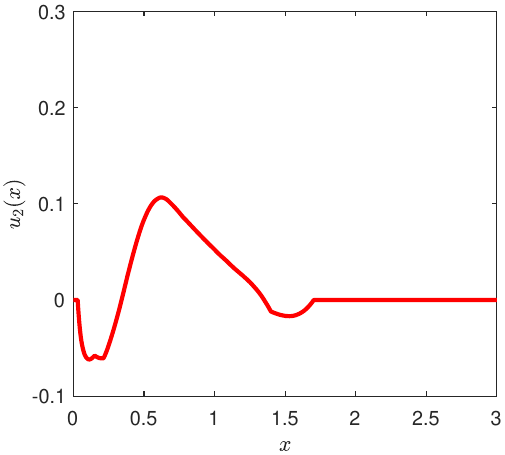}
\par\end{centering}
\begin{centering}
\includegraphics[scale=0.4]{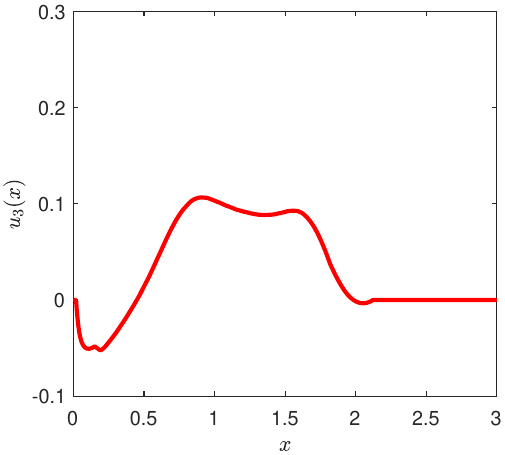}\includegraphics[scale=0.4]{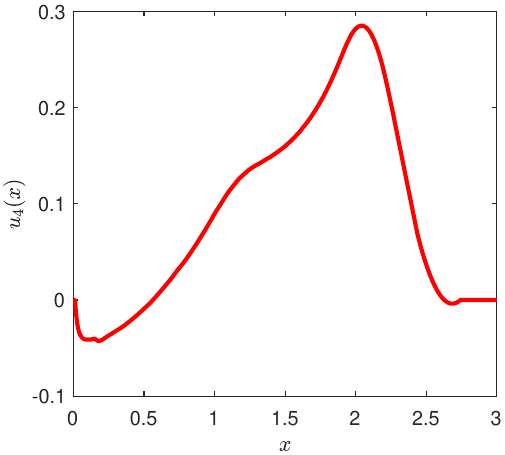}
\par\end{centering}
\caption{The optimal control sequence $u_{1},\cdots u_{4}$.\label{fig:control}}

\end{figure}

We simulate an example for the algorithm proposed in Section \ref{subsec:spec-case}.
For illustration, we consider a scalar system $x_{k+1}=x_{k}+0.3\sin x_{k}+u_{k}$
with Lagrangian $L(x,u)=0.01x^{4}+u^{2}$. Let the state space be
$X=[0,3]$ and the system evolves on $\mathcal{I}[1,5]$. The initial
and target distributions are $\rho_{1}\sim\mathscr{N}(0.7,0.03)$
and $\rho_{5}\sim\mathscr{N}(2.1,0.05)$ respectively. To implement
the algorithm, we discretize $X$ to uniform pieces. %
The result is shown in Fig. \ref{fig:cost} (optimal cost), \ref{fig:density}
(density evolution) and \ref{fig:control} (optimal controller). 

We can see from Fig. \ref{fig:cost} that the algorithm converges
after around $250$ iterations. In Fig. \ref{fig:density}, the two
shaded areas represent the initial and target densities respectively
and the rest of the curves represent the evolution of the densities
on $\mathcal{I}[2,5]$. We can see that $\rho_{5}$ (in purple color)
is close to the target. Fig. \ref{fig:control} shows the control
actions at each step. Since the initial density $\rho_{1}$ is almost
concentrated on the interval $[0,1.5]$, the control action $u_{1}$
is also concentrated on this interval. In the steps that follow, the
control actions gradually move to the right half of $X$, this is
due to that the target distribution is mainly concentrated on $[1.5,3]$. 


\vspace{-2mm}
\section{Conclusion}
We have proposed several equivalent formulations of dDOT. The duality
formulation is well-suited for numerical implementation. However,
the numerical algorithm is presently only valid for systems with full
input: $x_{k+1}=f(x_{k})+u_{k}$. Extending it to general systems
is the central task of future work. 

\vspace{-2mm}
\section{Appendix}
\emph{$c$-transform}: given a function $c:X\times Y\to\mathbb{R}$
and two functions $\phi:X\to\mathbb{R}\cup\{-\infty\}$, $\varphi:Y\to\mathbb{R}\cup\{-\infty\}$,
the $c$-transform of $\phi$ and $\varphi$ are defined as
\begin{align*}
\phi^{c}(y) & =\inf_{x\in X}\{c(x,y)-\phi(x)\}\\
\varphi^{c}(x) & =\inf_{y\in X}\{c(x,y)-\varphi(y)\}
\end{align*}
respectively.

\emph{The operator $K$} in (\ref{op:K}). We show that $K$ is a
bounded linear operator with norm $2\sqrt{\mu_{X}}$. To see this,
we compute{\small{}
\begin{align*}
\|K\|_{\mathcal{L}(\mathscr{H},\mathscr{G})}^{2} & =\sup_{\|v\|_{\mathscr{H}}\le1}\|Kv\|_{\mathscr{G}}^{2}\\
 & =\sup_{\|v\|_{\mathscr{H}}\le1}\sum_{k=1}^{T-1}\int_{X\times X}(v_{k}(x)-v_{k+1}(y))^{2}dxdy\\
 & \le2\mu_{X}\sup_{\|v\|_{\mathscr{H}}\le1}\sum_{k=1}^{T-1}\int_{X}v_{k}(x)^{2}dx+v_{k+1}(x)^{2}dx\\
 & =4\mu_{X}.
\end{align*}
}Now let $v_{1}=v_{T}=0$ and $v_{k}(x)=\frac{(-1)^{k}}{\sqrt{(T-2)\mu_{X}}},\;k\in\mathcal{I}[2,T-1]$.
Then obviously $\|v\|=1$, and thus {\small{}
\begin{align*}
\|K\|_{\mathcal{L}(\mathscr{H},\mathscr{G})}^{2} & \ge\sum_{k=2}^{T-1}\int_{X\times X}(v_{k}(x)-v_{k+1}(y))^{2}dxdy\\
 & =\sum_{k=2}^{T-1}\int_{X\times X}\frac{4}{(T-2)\mu_{X}}dxdy\\
 & =4\mu_{X}.
\end{align*}
}as desired. 
\begin{lem}
\label{lem:K-Gauss}For the Kantorovich problem 
\[
\mathcal{K}(\mu,\nu)=\min_{\pi\in\Gamma(\mu,\nu)}\int_{\mathbb{R}^{n}\times\mathbb{R}^{n}}c(x,y){\rm d}\pi(x,y),
\]
if the cost function $c(x,y)$ is quadratic in $[x^{\top},y^{\top}]^{\top}$,
and $\mu$, $\nu$ are zero mean Gaussian distributions, then 
\[
\mathcal{K}(\mu,\nu)=\max_{\phi,\varphi}\int_{\mathbb{R}^{n}}\phi{\rm d}\mu-\varphi{\rm d}\nu
\]
where the functions $\phi,\varphi$ are quadratic functions of the
form $x^{\top}Px$ for some symmetric matrices $P$. %

\end{lem}
\begin{IEEEproof}
Since $c(x,y)$ is quadratic, $\mathcal{K}(\mu,\nu)$ depends only
on the covariance matrix of $\pi$. Therefore, it is sufficient to
restrict $\pi$ to Gaussian distributions. Partition ${\rm Cov}(\pi)$
as
{\small
\[
{\rm Cov}(\pi)=\begin{bmatrix}M_{1} & M_{2}\\
M_{2}^{\top} & M_{3}
\end{bmatrix},
\]
}
we then have 
{\small
\begin{align*}
\int_{\mathbb{R}^{n}\times\mathbb{R}^{n}} & c(x,y){\rm d}\pi(x,y)=\int_{\mathbb{R}^{n}}x^{\top}M_{1}x{\rm d}\mu(x)\\
 & \quad+\int_{\mathbb{R}^{n}}y^{\top}M_{3}y{\rm d}\nu(y)+2\int_{\mathbb{R}^{n}\times\mathbb{R}^{n}}x^{\top}M_{2}y{\rm d}\pi(x,y).
\end{align*}
}
Let $\rho(x,y)$, $\rho(x)$, and $\rho(y|x)$ be the densities of
$\pi$, $\mu$ and the conditional distrubition $Y|X=x$ respectively.
Then
{\small
\begin{align*}
\int_{\mathbb{R}^{n}\times\mathbb{R}^{n}}x^{\top}M_{2}y{\rm d}\pi(x,y) & =\int_{\mathbb{R}^{n}}\int_{\mathbb{R}^{n}}x^{\top}M_{2}y\rho(y|x)\rho(x){\rm d}y{\rm d}x\\
 & =\int_{\mathbb{R}^{n}}x^{\top}M_{2}\mathbb{E}[Y|X=x]\rho(x){\rm d}x.
\end{align*}
}
Now since the conditional distribution $Y|X=x$ is Gaussian whose
mean value is of the form $Ux+v$ for some constant matrix $U$ and
vector $v$, the integral above must be quadratic in $x$ (remember
that $\mathbb{E}[\mu]=0$). Therefore, $\int_{\mathbb{R}^{n}\times\mathbb{R}^{n}}c{\rm d}\pi$
can be written as $\int_{\mathbb{R}^{n}}x^{\top}\tilde{M}_{1}x{\rm d}\mu(x)+\int_{\mathbb{R}^{n}}y^{\top}\tilde{M}_{3}y{\rm d}\nu(y)$.
At the minimum, there must hold $x^{\top}\tilde{M}_{1}x+y^{\top}\tilde{M}_{3}y\le c(x,y)$
invoking Kantorovich duality. 
\end{IEEEproof}
\begin{lem}
\label{lem:SDP}Let $A=A^{\top}\in\mathbb{R}^{n\times n}$, $c\in\mathbb{R}^{n}$,
$d\in\mathbb{R}$, then 
\begin{equation}
\min_{x}\{x^{\top}Ax+2c^{\top}x+d\}\ge0\label{eq:lem:min}
\end{equation}
is equivalent to
\begin{equation}
\begin{bmatrix}A & c\\
c^{\top} & d
\end{bmatrix}\ge0.\label{eq:lem:SDP}
\end{equation}
\end{lem}

\bibliographystyle{IEEEtran}
\bibliography{IEEEabrv,OT}

\end{document}